\documentclass[12pt]{article}
\usepackage{amssymb, amsmath}
\usepackage{color}
\begin{document}
\renewcommand{\theequation}{\arabic{section}.\arabic{equation}}
\newtheorem{thm}{Theorem}[section]
\newtheorem{lem}{Lemma}[section]
\newtheorem{deff}{Definition}[section]
\newtheorem{pro}{Proposition}[section]
\newtheorem{cor}{Corollary}[section]
\newcommand{\n}{\nonumber}
\newcommand{\vare}{\varepsilon}
\newcommand{\bb}{\begin{equation}}
\newcommand{\ee}{\end{equation}}
\newcommand{\bq}{\begin{eqnarray}}
\newcommand{\eq}{\end{eqnarray}}
\newcommand{\bqn}{\begin{eqnarray*}}
\newcommand{\eqn}{\end{eqnarray*}}

\title{On the temporal decay for the Hall-magnetohydrodynamic equations}
\author{Dongho Chae$^{(a)}$ and Maria Schonbek$^{(b)}$\\
\ \\
(a) Department of Mathematics\\
   Chung-Ang University\\
   Seoul156-756, Korea\\
   e-mail: dchae@cau.ac.kr\\
   \ \\
   (b) Department of Mathematics\\
   University of California\\
   Santa Cruz, CA 95064, USA\\
      e-mail: schonbek@ucsc.edu}

 \date{}
\maketitle
\begin{abstract}
We establish  temporal decay estimates for  weak solutions to the
Hall-magnetohydrodynamic equations. With these estimates in hand  we obtain algebraic  time decay  for higher order Sobolev norms of small initial data solutions.
\end{abstract}

\section{Introduction}
 \setcounter{equation}{0}
We consider the  incompressible MHD-Hall equations in $\Bbb R^3$.
\bq
& & \hspace{-1cm}  \partial_t u +  u \cdot \nabla u + \nabla p  =  (\nabla \times B) \times B + \nu\Delta u,  \label{u} \\
& & \hspace{-1cm} \nabla \cdot u = 0 , \quad \nabla \cdot B = 0,\label{eq:divu}\\
& & \hspace{-1cm} \partial_t B -  \nabla \times (u \times B) + \nabla \times ( (\nabla \times B) \times B)  = \mu\Delta B ,\label{B1}\\
& & \hspace{-1cm} u(x,0)=u_0(x)\quad;\quad B(x,0)=B_0(x). \label{B}
\eq Here  $u=(u_1,u_2,u,u_3)=u(x,t)$ is the velocity of the charged
fluid, $B=(B_1. B_2, B_3)$  the magnetic field induced by the motion
of the charged fluid, $p=p(x,t)$  the pressure of the fluid. The
positive constants $\nu$ and $\mu$ are the viscosity and the
resistivity coefficients. Without loss of generality  we let
$\mu=\nu=1$. Compared with the usual viscous incompressible MHD
system, the system (\ref{u})-(\ref{B}) contains the extra term
$\nabla \times ( (\nabla \times B) \times B)=\nabla \times ((\nabla
\times B)\times B)$, which is the so called  Hall term. This term is
important when the magnetic shear is large, where the magnetic
reconnection happens. On the other hand, in the case of laminar
flows where the shear is  weak, one ignores the Hall term, and the
system reduces to the usual MHD.  We refer \cite{for,sim} for the physical background of the magnetic reconnection and the Hall-MHD.\\
Compared to the case of the usual
MHD the history of the fully rigorous mathematical study of the
Cauchy problem for the Hall-MHD system is very short. The global
existence of weak solutions in the periodic domain is done in
\cite{ach} by  a Galerkin approximation. The global existence in
the whole domain in $\Bbb R^3$ as well as the local well-posedness
of smooth solution is proved in \cite{cha}, where the global
existence of smooth solution for small initial data is also
established.

  In this paper we  study the Cauchy problem of
the Hall-MHD system and establish  temporal decay estimates for
the solutions. Our results,  provide a mathematically rigorous
basis  to explain the decay of energy in the Hall-MHD, which had been obtained by numerical simulations (see e.g. \cite{chr,miu}).

Algebraic rates for the asymptotic behavior of solutions to the
Navier-Stokes equations were obtained first by the second author of
this paper in \cite{sch1}, using the method of Fourier Splitting. This technique  was introduced
first to study the decay of solutions to parabolic conservation laws \cite{sch0}. The Fourier Splitting method was then refined in in \cite{sch2, sch3}. (see also \cite{kaj}).
  Here we apply the arguments of \cite{sch1} and
\cite{sch2} to obtain algebraic time decay rates for the solutions
of the Hall-MHD system. The existence of Hall term in the equations
generates extra terms to control, which needed to be   handled
in our proofs  by introducing new estimates.  We now list the main theorems of the paper. The first  is a preliminary decay estimate for the
weak solutions.

\begin{thm}
Let $(u_0, B_0) \in (L^2(\Bbb R^3)\cap L^1(\Bbb R^3))^2 $ with div
$B_0=$ div $u_0=0$. Then, there exists a weak solution to the system
(\ref{u})-(\ref{B}), which satisfies \bb\label{th1} \|u(t)\|_{L^2}^2
+ \|B(t)\|_{L^2}^2\leq C(t+1)^{-\frac32}. \ee
\end{thm}

\noindent The following is a decay estimate for the higher order Sobolev
norms, whose global in time existence is guaranteed for sufficiently
small initial data(\cite{cha}).

\begin{thm}
Let $(u_0, B_0) \in ( L^1(\Bbb R^3))^2 $ satisfies the conditions of
Theorem 2.1 below.  Assume that
 \bb\label{th1a}
  \|u(t)\|_{L^2}^2 +\|B(t)\|_{L^2} ^2 \leq C_0 (t+1)^{-2\mu}
  \ee
  for $t\geq 0$ and $\mu \geq0$. Then, for $m\in \Bbb N$ there
  exists $C_m=C_m(\mu, C_0)$ and $T_* >0$ such that
  \bb\label{th14b}
  \|D^m u(t)\|_{L^2}^2 +\|D^m B(t)\|_{L^2} ^2 \leq C_m (t+1)^{-m-2\mu}
  \ee
for all $t\geq T_*$.
\end{thm}

\noindent{\em Remark 1.1 } Since (\ref{th1a}) is valid for $\mu=3/4$
by Theorem 1.1, we obtain the
 decay estimate,
 \bb\label{th14c}
  \|D^m u(t)\|_{L^2}^2 +\|D^m B(t)\|_{L^2} ^2 \leq C_m (t+1)^{-m-\frac32}
  \ee
  for $t\geq T_*$.
\section{Proof of the main theorems}

We recall that  to use the Fourier Splitting technique  we need the
following two  main estimates: Let  $V(\cdot,t)\in L^2(\Bbb R^n)\cap
L^1(\Bbb R^n)$,
\begin{align} \frac{d}{dt} \|V(t)\|_{L^2} \leq -C \|\nabla V(t)\|_{L^2}^2.\label{F1}\\
 |\widehat{V}(\xi,t) \leq C \;\;\;\mbox{for}\, |\xi|\ll 1,\label{F2}
\end{align}
where $ \widehat{V}$ denotes the Fourier transform of $V$ defined by
$$\widehat{V}(\xi)=\frac{1}{(2\pi)^{\frac{n}{2}}}\int_{\Bbb R^n} V(x)
e^{-ix\cdot\xi } dx, \qquad i=\sqrt{-1}.
$$
These conditions will insure that $\|V(t)\|_{L^2}$ decays at the
same rate as the solutions of the heat equations, with the same data
$V(x,0)$, \cite{schO}. We will use this method with appropriate
modifications for our equations. The next Lemma is in the spirit of
condition (\ref{F2}) above.
 \setcounter{equation}{0}
 \begin{lem}
 Let $(u, B)$ be a smooth solution to the system  (\ref{u})-(\ref{B}) with $(u_0, B_0)$ satisfying the initial condition as in Theorem 1.1. Then, we have
 \bb\label{lem1}
 |\hat{u}(\xi,t)| +|\hat{B}(\xi,t)|\leq C\left(1+\frac{1}{|\xi|} \right),
 \ee
 where $C=C(\|u_0\|_{L^1\cap L^2} +\|B_0\|_{L^1\cap L^2} )$.
 \end{lem}
 {\bf Proof } Using the elementary vector calculus, one can rewrite (\ref{u}) and (\ref{B1}) as
 \bb
 u_t-\Delta u= -\Bbb P \nabla \cdot(u\otimes u-B\otimes B),
 \ee
 and
 \bb
 B_t -\Delta B= - \nabla \cdot(u\otimes B-B\otimes u)-\nabla \times \,\{\nabla \cdot\, (B\otimes B)\}
  \ee
  respectively, where $\Bbb P$ is the Leray projection operator defined by $\Bbb P f=f-\nabla \Delta ^{-1} \nabla \cdot f$.
Hence, we have the following representation of solutions in terms of the Fourier transform,
 \bq
 &&\hat{u}(\xi,t)=e^{-|\xi|^2 t}\hat{u}_0(\xi)\n \\
  &&\quad -\int_0 ^t e^{-|\xi|^2(t-s)}\left(1-\frac{\xi\otimes\xi}{|\xi|^2}\right) \left\{ \xi\cdot (\widehat{u\otimes u}) (\xi,s) -\xi\cdot (\widehat{B\otimes B })(\xi,s) \right\}ds,\n \\
 \eq
 and
 \bq
 &&\hat{B}(\xi,t)=e^{-|\xi|^2 t}\hat{B}_0(\xi)\n \\
  &&\quad -\int_0 ^t e^{-|\xi|^2(t-s)} \left\{ \xi\cdot (\widehat{u\otimes B}) (\xi,s) -\xi\cdot (\widehat{B\otimes u })(\xi,s)-\xi \times \{ \xi\cdot (\widehat{B\otimes B })(\xi,s) \right\}ds.\n \\
 \eq
 From these representations we obtain
 \bq
 |\hat{u}(\xi,t)|&\leq& |\hat{u}_0 (\xi)|+ \int_0 ^t |\xi|  e^{-|\xi|^2(t-s)}\{ | (\widehat{u\otimes u}) (\xi,s)| +|(\widehat{B\otimes B })(\xi,s)|\} ds\n \\
 &\leq& \|u_0\|_{L^1} +\int_0 ^t |\xi|  e^{-|\xi|^2(t-s)}(\|u(s)\|_{L^2} ^2 +\|B(s)\|_{L^2}^2)ds\n \\
  &\leq& \|u_0\|_{L^1} + (\|u_0\|_{L^2} ^2 +\|B_0\|_{L^2}^2)\frac{1}{|\xi|} (1-e^{-|\xi|^2t})\n \\
  &\leq& C \left(1+\frac{1}{|\xi|}\right),
  \eq
  and
  \bq
 |\hat{B}(\xi,t)|&\leq& |\hat{B}_0 (\xi)|+ \int_0 ^t  e^{-|\xi|^2(t-s)}|\xi|\{ | (\widehat{u\otimes B}) (\xi,s)| +|(\widehat{B\otimes u })(\xi,s)|\}ds\n \\
 &&\quad    + \int_0 ^t  e^{-|\xi|^2(t-s)}|\xi|^2|(\widehat{B\otimes B })(\xi,s)| \} ds\n \\
 &\leq& \|B_0\|_{L^1} +\int_0 ^t  e^{-|\xi|^2(t-s)}(2|\xi|\|u(s)B(s)\|_{L^1}  +|\xi|^2\|B(s)\|_{L^2}^2)ds\n \\
  &\leq& \|B_0\|_{L^1} + \left\{\frac{2}{|\xi|}\|u_0\|_{L^2}\|B_0\|_{L^2}+ \|B_0\|_{L^2}^2)\right\} (1-e^{-|\xi|^2t})\n \\
  &\leq& C \left(1+\frac{1}{|\xi|}\right).
  \eq
  $\square$\\
  \ \\
  \noindent{\bf Proof of Theorem 1.1 }  Multiplying (\ref{u}) by $u$, and (\ref{B1}) by $B$,
  and integrating over $\Bbb R^3$, and integrating by part we get
  \bb\label{energy}
  \frac{d}{dt} \int_{\Bbb R^3} ( |u|^2 +|B|^2 )dx=-2\int_{\Bbb R^3} (|\nabla u|^2 +|\nabla B|^2)dx.
  \ee
  Using the estimates (\ref{energy}) and (\ref{lem1}), and applying the standard Fourier-Splitting method, developed in \cite{sch1} and \cite{sch2},  one can conclude the estimate
  (\ref{th1}). More specifically  by (\ref{energy}) and (\ref{lem1}), one obtains a preliminary decay estimate that allows, as was done for the Navier-Stokes equations, to obtain the better estimate on the Fourier transfer of $(u,B)$ near the origin in frequency space,  as required in condition (\ref{F2}).\;\;\;$\square$

  \bigskip

  In order to prove Theorem 1.2
we first recall the following small data global regularity result
proved in \cite{cha}(see Theorem 2.3).
\begin{thm}
Let $m\in \Bbb N, m\geq 3$ and $(u_0, B_0) \in [H^m (\Bbb R^3)]^2$
with div $B_0=$ div $u_0=0$. There exists a constant $K_1$ such that
if $\|u_0\|_{H^m} +\|B_0\|_{H^m} \leq K_1$, then there exists a
unique solution $(u,B) \in L^\infty (\Bbb R_+; H^m (\Bbb R^3))$
satisfying \bb\label{cdl} \frac{d}{dt} (\|u\|_{H^m}^2 +\|B\|_{H^m}^2
) \leq -(\|\nabla u\|_{H^m}^2 +\|\nabla B\|_{H^m}^2 ). \ee
\end{thm}

\noindent{\em Remark 2.1 } Although the global energy inequality
(\ref{cdl}) is not written down in \cite{cha}, it is immediate by
choosing the constant $K_1=\frac12 K$, where $K$ is the constant
in Theorem 2.3(\cite{cha}) bounding the initial data to obtain the global smooth solution.\\
\ \\
We observe the following fact.
  \begin{lem}
  Let $(u_0, B_0)\in ( L^1(\Bbb R^3))^2$ be as in Theorem 2.1.
  Then, for all $|\xi|\leq 1$ and for all $j\in \Bbb N$ we have
  \bb
  |\widehat{D^j u}(\xi,t)|+  |\widehat{D^j B}(\xi ,t)|\leq C,
  \ee
 where $C=C(\|u_0\|_{L^1\cap L^2} +\|B_0\|_{L^1\cap L^2} )$.
  \end{lem}
  {\bf Proof } Let $|\xi|\leq 1$.  Using the result of Lemma 2.1, we have
  \bqn
   \lefteqn{|\widehat{D^j u}(\xi,t)|+|\widehat{D^j B}(\xi ,t)|
   \leq |\xi|^j (|\hat{u}(\xi,t)|+|\hat{B}(\xi,t)|)}\n \\
   &&\qquad \leq |\xi|(|\hat{u}(\xi,t)|+|\hat{B}(\xi,t)|)\leq C |\xi| \left(1+\frac{1}{|\xi|} \right)
   \leq C.
   \eqn
  $\square$\\
  \ \\
  The following is an auxiliary decay estimate for higher order Sobolev norms.
  \begin{thm}
 Let $(u_0, B_0)\in ( L^1(\Bbb R^3))^2$ be as in Theorem 2.1. Then, there exists a constant $C$ such that
 \bb\label{prop1}
 \|u(t)\|_{H^m}^2 +\|B(t)\|_{H^m}^2 \leq C(t+1)^{-\frac32}.
 \ee
\end{thm}
{\bf Proof } We apply the Fourier-Splitting method. Let $(\alpha_1,
\alpha_2, \alpha_3)\in [\Bbb N \cup\{0\}]^3$,
$|\alpha|=\alpha_1+\alpha_2+\alpha_3$, be the multi-index. The
Fourier transform of (\ref{cdl}) is written as \bq\label{14}
 \lefteqn{\frac{d}{dt} \int_{\Bbb R^3} \left(\sum_{|\alpha|\leq m }|\widehat{D^\alpha
 u}|^2 +\sum_{|\alpha|\leq m }|\widehat{D^\alpha
 B}|^2\right) d\xi}\n \\
 && = \frac{d}{dt} \int_{\Bbb R^3} \sum_{|\alpha|\leq m }
 (|\xi_1|^{\alpha_1}  |\xi_2|^{\alpha_2} |\xi_3|^{\alpha_3})^2 ( |
 \hat{u}|^2 + |\hat{B} |^2) d\xi\n\\
 &&\qquad -\int_{\Bbb R^3} \sum_{|\alpha|\leq m }
 |\xi|^2(|\xi_1|^{\alpha_1}  |\xi_2|^{\alpha_2} |\xi_3|^{\alpha_3})^2 ( |
 \hat{u}|^2 + |\hat{B} |^2) d\xi.
 \eq
 Let
 $$S:=\left\{ \xi \in \Bbb R^3 \, |\, |\xi|\leq \left(
 \frac{k}{t+1}\right)^{\frac12}\right\}.
 $$
 From (\ref{14}) it follows that
 \bq
 \lefteqn{\frac{d}{dt} \int_{\Bbb R^3} (\sum_{|\alpha|\leq m }|\widehat{D^\alpha
 u}|^2 +\sum_{|\alpha|\leq m }|\widehat{D^\alpha
 B}|^2) d\xi}\n \\
&&\leq -\int_{S} |\xi|^2(\sum_{|\alpha|\leq m }|\widehat{D^\alpha
 u}|^2 +\sum_{|\alpha|\leq m }|\widehat{D^\alpha
 B}|^2) d\xi \n \\
 &&\qquad -\int_{S^c} |\xi|^2(\sum_{|\alpha|\leq m
}|\widehat{D^\alpha
 u}|^2 +\sum_{|\alpha|\leq m }|\widehat{D^\alpha
 B}|^2t) d\xi\n \\
 &&\leq -\frac{k}{t+1}\int_{\Bbb R^3}(\sum_{|\alpha|\leq m
}|\widehat{D^\alpha
 u}|^2 +\sum_{|\alpha|\leq m }|\widehat{D^\alpha
 B}|^2) d\xi\n \\
 &&\qquad+\frac{k}{t+1}\int_{S} (\sum_{|\alpha|\leq m
}|\widehat{D^\alpha
 u}|^2 +\sum_{|\alpha|\leq m }|\widehat{D^\alpha
 B}|^2) d\xi.\n \\
 \eq
 Since $|\xi|\leq k (t+1)^{-\frac12}$, and $|\widehat{D^\alpha
 u}(\xi,t)| +|\widehat{D^\alpha
 B}(\xi,t)|\leq C $ for $\xi \in S$ and $t\geq T_0$ for some
 $T_0>0$ by Lemma 2.2, we have
\bq
  \lefteqn{\frac{d}{dt} \left[(1+t)^k\int_{\Bbb R^3}(\sum_{|\alpha|\leq m }|\widehat{D^\alpha
 u}|^2 +\sum_{|\alpha|\leq m }|\widehat{D^\alpha
 B}|^2) d\xi\right] }\n \\
 &&\quad\leq k (t+1)^{k-1}\int_S(\sum_{|\alpha|\leq m
}|\widehat{D^\alpha
 u}|^2 +\sum_{|\alpha|\leq m }|\widehat{D^\alpha
 B}|^2) d\xi\n \\
 &&\quad\leq  C(t+1)^{k-\frac52} \qquad \forall t\geq T_0.
 \eq
 Integrating over $[T_0, t]$, and dividing by $(t+1)^k$,
 we find
\bq
 &&\int_{\Bbb R^3}(\sum_{|\alpha|\leq m }|\widehat{D^\alpha
 u}(\xi,t)|^2 +\sum_{|\alpha|\leq m }|\widehat{D^\alpha
 B}(\xi,t)|^2) d\xi\n \\\
 &&\qquad\leq (t+1)^{-k}\int_{\Bbb R^3}(\sum_{|\alpha|\leq m }|\widehat{D^\alpha
 u}(\xi,T_0)|^2 +\sum_{|\alpha|\leq m }|\widehat{D^\alpha
 B}(\xi, T_0)|^2) d\xi  \n \\
 &&\quad\qquad +C (t+1)^{-\frac32}.
 \eq
 The estimate (\ref{prop1}) follows if we choose $k=3/2$.
 $\square$\\
 \ \\
 In order to establish Theorem 1.4 we first show the following
 auxiliary lemma, which is similar in form to Lemma 3.2.
\begin{lem}
Let $(u, B)$ is a smooth solution of (\ref{u})-(\ref{B}) and $m\in
\Bbb N$. Then, we have the following inequality.
 \bq\label{lem23}
  &&\frac{d}{dt} (\|D^m u\|_{L^2}^2 + \|D^m B\|_{L^2}^2) + \|D^{m+1} u\|_{L^2}^2 + \|D^{m+1}
  B\|_{L^2}^2\n \\
  &&\quad\leq C_m
  (\|u\|_{L^\infty}^2 \|D^m u\|_{L^2}^2 +\|B\|_{L^\infty} ^2 \|D^m
  B\|_{L^2} ^2+\|B\|_{L^\infty} ^2 \|D^m u\|_{L^2}^2\n \\
  &&\qquad\qquad\qquad
  +\|u\|_{L^\infty} ^2 \|D^m B\|_{L^2} ^2 +\|\nabla B
  \|_{L^\infty}^2
  \|D^m B\|_{L^2} ^2 )+R_m, \n \\
  \eq
 where
 \bb\label{rm}
 R_m=\left\{ \aligned &\qquad0, \qquad\quad\mathrm{if}\qquad m=1,2\\
  &C_m\sum_{1\leq j\leq m/2} (\|D^j u\|_{L^\infty} ^2 \|D^{m-j} u\|_{L^2}
  ^2 +\|D^j B\|_{L^\infty} ^2 \|D^{m-j} B\|_{L^2}^2)\\
   &\quad+C_m\sum_{1\leq j\leq m/2} (\|D^j B\|_{L^\infty} ^2 \|D^{m-j} u\|_{L^2}^2
  +\|D^j u\|_{L^\infty} ^2 \|D^{m-j} B\|_{L^2}^2)\\
  &\quad+C_m\sum_{2\leq i\leq \frac{m+1}{2}} \|D^i B\|_{L^\infty}^2 \|D^{m+1-i} B\|_{L^2} ^2,
   \quad\mathrm{if}\qquad m\geq 3.
  \endaligned
  \right.
  \ee
\end{lem}
{\em Remark 2.2 } Although the Hall term generates the extra factor
of norm for the derivative of $B$ such as $\|\nabla B\|_{L^\infty}$,
which is not present in \cite{sch3}, this can be handled without
difficulty as shown in the proof of Theorem 1.2 below.\\
\ \\
 \noindent{\bf Proof of Lemma 2.3  } Let $m\geq 3$.  Multiplying (\ref{u})
and (\ref{B1}) by $u$ and $B$ respectively, and integrating each one
over $\Bbb R^3$, and integrating by part, we obtain the following
inequalities.
 \bq\label{lem23a}
 \frac{d}{dt} \int_{\Bbb R^3} |D^m u |^2 dx&\leq& C\sum_{j,k=1}^3
 \int_{\Bbb R^3} |D^m (u_k u_j)|^2 dx +C\sum_{j,k=1}^3
 \int_{\Bbb R^3} |D^m (B_k B_j)|^2 dx\n \\
 &&\quad-
 \int_{\Bbb R^3} |D^{m+1} u|^2 dx,
 \eq
 and
\bqn
 \frac{d}{dt} \int_{\Bbb R^3} |D^m B |^2 dx&\leq& C\sum_{j,k=1}^3
 \int_{\Bbb R^3} |D^m (B_k u_j)|^2 dx +C\sum_{j,k=1}^3
 \int_{\Bbb R^3} |D^m (u_k B_j)|^2 dx\n \\
 &&\quad+C \int_{\Bbb R^3} |D^m (B\cdot \nabla) B |^2 dx-
 \int_{\Bbb R^3} |D^{m+1} B|^2 dx,
 \eqn
 where we used the fact $\|\nabla f\|_{L^2}= \|\nabla \times f\|_{L^2}$ if $\nabla \cdot f=0$.
We have the following auxiliary estimates. For $m\geq 2$ \bq
 \sum_{j,k=1}^3\int_{\Bbb R^3} |D^m (u_ku_j )|^2dx&\leq& C
\sum_{j,k=1}^3\sum_{1\leq i\leq m/2} \int_{\Bbb R^3} |D^{m-i} u_k D^i u_j|^2 dx+C \int_{\Bbb R^3} |D^m u|^2 |u|^2 dx\n \\
&\leq& C \sum_{1\leq i\leq m/2} \|D^i u\|_{L^\infty}^2 \|D^{m-i} u\|_{L^2}^2 + C \|u\|_{L^\infty}^2 \|D^m u\|_{L^2}^2.\n \\
\eq
Similarly, we have
\bb
 \sum_{j,k=1}^3\int_{\Bbb R^3} |D^m (B_k B_j )|^2dx\leq C \sum_{1\leq i\leq m/2} \|D^i B\|_{L^\infty}^2 \|D^{m-i} B\|_{L^2}^2 + C \|B\|_{L^\infty}^2 \|D^m B\|_{L^2}^2,
  \ee
 and
 \bq
 &&\sum_{j,k=1}^3\int_{\Bbb R^3} |D^m (B_k u_j )|^2dx+\sum_{j,k=1}^3\int_{\Bbb R^3} |D^m (u_k B_j )|^2dx\n \\
 &&\quad\leq C \sum_{1\leq i\leq m/2} \|D^i B\|_{L^\infty}^2 \|D^{m-i} u\|_{L^2}^2 + C \|B\|_{L^\infty}^2 \|D^m u\|_{L^2}^2\n \\
&&\qquad + C \sum_{1\leq i\leq m/2} \|D^i u\|_{L^\infty}^2 \|D^{m-i} B\|_{L^2}^2 + C \|u\|_{L^\infty}^2 \|D^m B\|_{L^2}^2.
  \eq
 For the Hall term we have the estimate
 \bq
  \int_{\Bbb R^3} |D^m (B\cdot \nabla B)|^2 dx &=&  \int_{\Bbb R^3} |D^m \nabla \cdot (B\otimes B)|^2 dx \leq \sum_{j,k=1}^3\int_{\Bbb R^3} |D^{m+1}(B_j B_k)|^2 dx\n \\
  &\leq& \|D^{m+1} B\|_{L^2}^2 \|B\|_{L^\infty}^2
   +\|D^m B\|_{L^2} ^2 \|\nabla B\|_{L^\infty}^2\n \\
   &&\quad +\sum_{2\leq i\leq \frac{m+1}{2}}  \|D^i B\|^2_{L^\infty} \|D^{m+1-i} B\|_{L^m} .
   \eq
  Let $k\geq 3$, By  the Gagliardo-Nirenberg inequality and Theorem 1.1 and Theorem 1.3 we have
  \bb
  \label{gag1}\|B\|_{L^\infty} \leq C\|B\|_{L^2} ^{\frac{2k-3}{2k}} \|D^k B\|_{L^2} ^{\frac{3}{2k}} \leq C(t+1)^{-\frac34}
   \quad \forall t\geq T_1,
  \ee
  for some $T_1 >0$, where $\theta\in (0, 1)$.  Hence, for any there
  exists $T_1>0$ such that
  $\vare>0$
  \bb\label{lem23b}
  \|D^{m+1} B(t)\|_{L^2}^2 \|B(t)\|_{L^\infty}^2 \leq \vare \|D^{m+1} B(t)\|_{L^2}^2.
 \ee
  for all $t\geq T_1$. Hence, this term can be absorbed into the viscosity term.
  Combining (\ref{lem23a})-(\ref{lem23b}) yields for $t\geq \max\{T_0, T_1\}$
  \bq
  &&\frac{d}{dt} (\|D^m u\|_{L^2}^2 + \|D^m B\|_{L^2}^2) +\|D^{m+1}u\|_{L^2}^2 +\|D^{m+1} B\|_{L^2}^2 \n \\
  &&\leq C_m
  (\|u\|_{L^\infty}^2 \|D^m u\|_{L^2}^2 +\|B\|_{L^\infty} ^2 \|D^m
  B\|_{L^2} ^2+\|B\|_{L^\infty} ^2 \|D^m u\|_{L^2}^2\n \\
  &&\qquad
  \qquad+\|u\|_{L^\infty} ^2 \|D^m B\|_{L^2} ^2 +\|\nabla B
  \|_{L^\infty}^2
  \|D^m B\|_{L^2} ^2 ) \n \\
  &&\qquad+C_m\sum_{1\leq j\leq m/2} (\|D^j u\|_{L^\infty} ^2 \|D^{m-j} u\|_{L^2}
  ^2 +\|D^j B\|_{L^\infty} ^2 \|D^{m-j} B\|_{L^2}^2)\n\\
   &&\qquad+C_m\sum_{1\leq j\leq m/2} (\|D^j B\|_{L^\infty} ^2 \|D^{m-j} u\|_{L^2}^2
  +\|D^j u\|_{L^\infty} ^2 \|D^{m-j} B\|_{L^2}^2)\n \\
  &&\qquad+C_m \sum_{2\leq i\leq \frac{m+1}{2}} \|D^i B\|_{L^\infty}^2 \|D^{m+1-i} B\|_{L^2} ^2.
  \eq
  This completes the proof of the lemma for $m\geq 3$.
Next we consider the case $m=1,2$. The estimate for $m=1$
corresponding to $u\cdot \nabla u$ is the following.
  \bqn
 \sum_{i,j,k=1}^3\int_{\Bbb R^3} \partial_i (u_j\partial_j u_k )
  \partial_i u_k dx&=&\sum_{i,j,k=1}^3\int_{\Bbb R^3} u_j \partial_i\partial_j
  u_k\partial_j u_k dx +\sum_{i,j,k=1}^3\int_{\Bbb R^3} \partial_iu_j
  u_k \partial_iu_k dx \n \\
  &&:=I_1+I_2.
  \eqn
For  $\vare >0$ we have
  $$ |I_1 |\leq \int_{\Bbb R^3} |u||Du| |D^2u|dx\leq C_\vare
  \|u\|_{L^\infty}^2 \|Du\|_{L^2}^2 +\vare \|D^2 u\|^2_{L^2}.
  $$
  Integrating by part for $I_2$ we have similar estimate,
$$ |I_2 |\leq C_\vare
  \|u\|_{L^\infty}^2 \|Du\|_{L^2}^2 +\vare \|D^2 u\|^2_{L^2}.
  $$
 In the case  $m=2$ we have
\bqn
 \lefteqn{\sum_{i,j,k,m=1}^3\int_{\Bbb R^3} \partial_i \partial_j (u_k\partial_k u_m )
  D^2 u_m dx=\sum_{i,j,k,m=1}^3\int_{\Bbb R^3} u_k
  \partial_i\partial_j\partial_k u_m D^2 u_m dx}\n \\
  && +2\sum_{i,j,k,m=1}^3\int_{\Bbb R^3}
  \partial_i u_k \partial_j\partial_k u_m D^2 u_m dx
+\sum_{i,j,k,m=1}^3\int_{\Bbb R^3} \partial_i
   \partial_j u_k\partial_k u_m
  D^2 u_m dx \n \\
  &&:=J_1+J_2+J_3.
  \eqn
Let  $\vare >0$,  we have simple estimates
$$
 |J_1|\leq \int_{\Bbb R^3} |u| |D^2u| |D^3 u|dx \leq
 C_\vare\|u\|_{L^\infty} ^2 \|D^2 u\|_{L^2}^2
 +\vare\|D^3u\|_{L^2}^2.
 $$
For $J_2$ and $J_3$ we obtain similar estimates after moving, by
integration by part, the derivative of $\partial_i u_k$ and
$\partial_k u_m$ respectively to the other factors in the
integrands. Therefore we obtain
$$
 |J_2|+|J_3| \leq
 C_\vare\|u\|_{L^\infty} ^2 \|D^2 u\|_{L^2}^2
 +\vare\|D^3u\|_{L^2}^2.
 $$
 The estimates for the other terms corresponding to $B\cdot\nabla
 B, u\cdot\nabla  B, B\cdot\nabla  u$ are similar, and we omit them.
 We are left only with the Hall term.
 For $m=1$ we have
 \bqn
 \lefteqn{\sum_{i=1}^3\int_{\Bbb R^3} \partial_i \nabla \times (B\cdot
 \nabla B) \cdot \partial_i B dx=\sum_{i=1}^3\int_{\Bbb R^3}\partial_i(B\cdot
 \nabla B) \cdot \partial_i (\nabla \times B) dx}\n \\
 &&=\sum_{i=1}^3\int_{\Bbb R^3} (\partial_i B\cdot \nabla B) \cdot \partial_i (\nabla \times
 B)dx +\sum_{i=1}^3\int_{\Bbb R^3}  (B\cdot \nabla\partial_i B )\cdot \partial_i (\nabla \times
 B)dx\n \\
 &&\leq C_\vare\|\nabla B\|_{L^\infty}^2 \|\nabla B\|_{L^2}^2 +\vare \|D^2
 B\|_{L^2}^2 + \|B\|_{L^\infty} \|D^2 B\|_{L^2}^2\n \\
 &&\leq C_\vare\|\nabla B\|_{L^\infty}^2 \|\nabla B\|_{L^2}^2 +2\vare \|D^2
 B\|_{L^2}^2, \qquad \forall t\geq T_1
  \eqn
where we used the fact (\ref{gag1}),
 $\|B(t)\|_{L^\infty} \leq \vare $ for  all $t\geq T_1$. In the case
 $m=2$ we have
 \bqn
 \lefteqn{\sum_{i,j =1}^3\int_{\Bbb R^3} \partial_i \partial_j  \nabla \times (B\cdot
 \nabla B) \cdot D^2 B dx=\sum_{i,j=1}^3\int_{\Bbb R^3}\partial_i\partial_j(B\cdot
 \nabla B) \cdot D^2(\nabla \times B) dx}\n \\
 &&=\sum_{i,j =1}^3\int_{\Bbb R^3} (\partial_i \partial_j B\cdot \nabla B) \cdot D^2(\nabla \times
 B)dx +2\sum_{i,j =1}^3\int_{\Bbb R^3}  (\partial_i B\cdot \nabla\partial_j B )\cdot D^2(\nabla \times
 B)dx\n \\
 &&\qquad +\sum_{i,j =1}^3\int_{\Bbb R^3} ( B\cdot \nabla \partial_i \partial_j B) \cdot D^2(\nabla \times
 B)dx \n \\
  &&\leq C_\vare\|\nabla B\|_{L^\infty}^2 \|D^2 B\|_{L^2}^2 +\vare \|D^3
 B\|_{L^2}^2 + \|B\|_{L^\infty} \|D^3 B\|_{L^2}^2\n \\
 &&\leq C_\vare\|\nabla B\|_{L^\infty}^2 \|D^2 B\|_{L^2}^2 +2\vare \|D^3
 B\|_{L^2}^2 \qquad \forall t\geq T_1.
  \eqn
This completes the proof of the lemma for $m=1,2$. $\square$\\
\ \\
Next lemma is more or less a repetition of Lemma 3.3 of \cite{sch3}.
\begin{lem}
Let $m\in \Bbb N$, $t\geq T_*=\max\{ T_0, T_1\}$, where $T_0, T_1$
are the given in the above lemma. Assume
\bb
 \|D^{m-1} u\|_{L^2}^2 +\|D^{m-1} B\|_{L^2}^2\leq C_{m-1}
 (t+1)^{-\rho_{m-1}} \quad \forall t\geq T_*.
 \ee
 Suppose
 \bb\label{hyp2}
 \frac{d}{dt} \left( \|D^m u\|_{L^2}^2 +\|D^m B\|_{L^2} ^2\right)
 \leq C_0 (t+1)^{-1} \|D^m u\|_{L^2} ^2 +\sum_{i=1}^m C_i
 (t+1)^{-s_i}
 \ee
 with $s_i\geq \rho_{m-1} +2$. Then,
\bb
 \|D^{m} u\|_{L^2}^2 +\|D^{m} B\|_{L^2}^2\leq C_{m}
 (t+1)^{-\rho_{m}} \quad \forall t\geq T_*
 \ee
 with $\rho_m=1+\rho_{m-1}$, where $C_m=C_m (C_{m-1}, C_i, s_i,
 \rho_{m-1}, m)$ does not depend on $T_*$.
\end{lem}
{\bf Proof } We use the Fourier-Splitting argument. Let
$$ S=\left\{ \xi\in \Bbb R^3 \, |\, |\xi| \leq \left(\frac{ C_0
+k}{t+1}\right)^{\frac12} \right\}, \quad k=1+\max_{1\leq i\leq m}
\{ s_i\}.
$$
Then,
 \bqn
 \lefteqn{\|D^{m+1} u\|_{L^2}^2 +\|D^{m+1} B\|_{L^2}^2 \geq
 \int_{S^c} |\xi|^2 (|\widehat{D^m u}|^2 +|\widehat{D^m B}|^2 )d\xi }\n \\
 && \quad \geq \frac{C_0 +k}{t+1} (\|D^m u\|_{L^2}^2 +\|D^m
 B\|_{L^2}^2) -\frac{C_0 +k}{t+1}\int_S(|\widehat{D^m u}|^2 +|\widehat{D^m B}|^2
 )d\xi\n \\
&& \quad \geq \frac{C_0 +k}{t+1} (\|D^m u\|_{L^2}^2 +\|D^m
 B\|_{L^2}^2) -\left(\frac{C_0 +k}{t+1}\right)^2\int_S(|\widehat{D^{m-1} u}|^2 +|\widehat{D^{m-1} B}|^2
 )d\xi.
 \eqn
Using this last inequality and the hypothesis (\ref{hyp2}), we have
 \bqn
 \lefteqn{\frac{d}{dt} \left( \|D^m u\|_{L^2}^2 +\|D^m B\|_{L^2} ^2\right)
 +\frac{k}{t+1}( \|D^{m+1} u\|_{L^2}^2 +\|D^{m+1} B\|_{L^2}^2)}\n \\
 && \quad \leq C_{m-1} \frac{(C_0 +k)^2}{(t+1)^{2+\rho_{m-1}}}
 +\sum_{i=1}^m C_i (t+1)^{-s_i}.
 \eqn
 Multiplying $(t+1)^k$ and integrating in time, and dividing by
 $(t+1)^k$, we find
 $$
 \|D^m u\|_{L^2}^2 +\|D^m B\|_{L^2} ^2 \leq C_m
 (t+1)^{-(1+\rho_{m-1})} +\sum_{i=1}^m C_i (t+1)^{-s_i +1}.
  $$
  Since $s_i \geq \rho_{m-1} +2$, the conclusion of the lemma
  follows. $\square$\\
  \ \\
We are now ready to prove Theorem 1.2.\\
\ \\
\noindent{\bf Proof of Theorem 1.2 } We first consider the case
$m=1,2$. From the estimate (\ref{gag1}) for $B$ and its placements
by $u$ and $\nabla B$, and $k=2$ we have \bb\label{m12}
\|u(t)\|_{L^\infty}^2 +\|B(t)\|_{L^\infty}^2+\|\nabla
B(t)\|_{L^\infty}^2 \leq C(t+1)^{-\frac32} \ee
 for $t\geq T_0$. Substituting (\ref{m12}) into (\ref{lem23}), we obtain
 \bqn
  &&\frac{d}{dt} (\|D^m u\|_{L^2}^2 + \|D^m B\|_{L^2}^2) + \|D^{m+1} u\|_{L^2}^2 + \|D^{m+1}
  B\|_{L^2}^2\n \\
  &&\quad\leq C
  (t+1)^{-\frac32}( \|D^m u\|_{L^2}^2 + \|D^m
  B\|_{L^2} ^2 )
  \eqn
  for $m=1,2$. We can now apply Lemma 2.4 directly to obtain (\ref{th14b}).
 For $m\geq 3$ we need to estimate $R_m$ of (\ref{rm}).
 By the Gagliardo-Nirenberg inequality we have
 \bb\label{gag2}
 \|D^j u\|_{L^\infty} \leq C \|D^{m+1} u\|_{L^2} ^{a_j} \|u\|_{L^2} ^{1-a_j},
 \ee
 and
  \bb\label{gag3}
 \|D^j B\|_{L^\infty} \leq C \|D^{m+1} B\|_{L^2} ^{a_j} \|B\|_{L^2} ^{1-a_j},
 \ee
 where $ a_j=\frac{j+\frac32}{m+1}$.
  Substituting (\ref{gag2}) and (\ref{gag3}) into $R_m$ of (\ref{rm}), and using Young's inequality($ab\leq \frac{a^p}{p}+\frac{b^q}{q}, 1/p+1/q=1;a,b>0; p,q\geq 1$), we obtain
  \bq\label{rm1}
  \lefteqn{R_m \leq \frac18 (  \|D^{m+1} u\|_{L^2}^2 + \|D^{m+1}
  B\|_{L^2}^2)}\n \\
   &&\quad+C ( \|u\|_{L^2}^2 +\|B\|_{L^2}^2 )\left( \sum_{1\leq j\leq m/2} \|D^{m-j} u\|_{L^2} ^{\frac{2}{1-a_j}} +\sum_{1\leq j\leq m/2} \|D^{m-j} B\|_{L^2} ^{\frac{2}{1-a_j}}\right)\n \\
   && \leq \frac18 (  \|D^{m+1} u\|_{L^2}^2 + \|D^{m+1}
  B\|_{L^2}^2)+ C\sum_{1\leq j\leq m/2} (t+1)^{-s_j},
  \eq
  where $s_j = 2\mu + \frac{(m+1)(m-j)}{m-j-1/2}.$ Note that since $1\leq j\leq m/2$
  we have $s_j\geq 2\mu +m+1$. Substituting (\ref{gag2}), (\ref{gag3}) and (\ref{rm1}) into
  (\ref{lem23}),  we obtain the hypothesis (\ref{hyp2}).
  Applying Lemma 2.4 directly, we obtain the conclusion of the theorem for $m\geq 3$.
  $\square$\\
  \ \\
  $$ \mbox{\bf Acknowledgements } $$
  Part of this research is done while the first author
 was visiting UCSC. His research is supported partially by NRF Grant no.
 2006-0093854. The second author's research  is partially supported by the NSF Grant DMS-0900909.
 
\end{document}